\documentclass{gtart_a}
\pdfoutput=1

\title{Three-manifolds, virtual homology, and group determinants} 

\author{Daryl Cooper}
\givenname{Daryl}
\surname{Cooper}
\address{Math Department\\UCSB\\Santa Barbara, CA 93106\\USA}
\email{cooper@math.ucsb.edu}
\urladdr{}

\author{Genevieve S Walsh}
\givenname{Genevieve S}
\surname{Walsh}
\address{{\rm[GW]}\qua Department of Math\\Tufts University\\Medford, MA 02155
\\USA\\and\\\newline
D\'epartement de Math\'ematiques\\UQAM\\Montr\'eal, QC H3C 3J7\\Canada}
\email{genevieve.walsh@tufts.edu}
\urladdr{}

\volumenumber{10}
\issuenumber{}
\publicationyear{2006}
\papernumber{50}
\lognumber{0722}
\startpage{2247}
\endpage{2269}

\doi{}
\MR{}
\Zbl{}

\arxivreference{math.GT/0603152}


\received{8 March 2006}
\revised{16 October 2006}
\accepted{24 October 2006}
\published{29 November 2006}
\publishedonline{29 November 2006}
\proposed{David Gabai}
\seconded{Cameron Gordon, Joan Birman}
\corresponding{}
\editor{CPR}
\version{}

\AtBeginDocument{\let\bar\wbar\let\tilde\wwtilde\let\hat\what}
\makeop{nullity}\makeop{incl}\makeop{Ker}\makeop{stab}\makeop{sym}
\makeop{Im}\makeop{rank}\makeop{trace}\makeop{End}\makeop{Aut}

\makeatletter
\def\cnewtheorem#1[#2]#3{\newtheorem{#1}{#3}[section]
\expandafter\let\csname c@#1\endcsname\c@theorem}

\newtheorem{theorem}{Theorem}[section]
\cnewtheorem{lemma}[theorem]{Lemma}
\cnewtheorem{corollary}[theorem]{Corollary} 
\cnewtheorem{definition}[theorem]{Definition} 
\cnewtheorem{proposition}[theorem]{Proposition}
\cnewtheorem{remark}[theorem]{Remark} 
\cnewtheorem{question}[theorem]{Question} 
\cnewtheorem{example}[theorem]{Example} 
\cnewtheorem{addendum}[theorem]{Addendum} 
\cnewtheorem{scholium}[theorem]{Scholium} 

\makeatother  

\let\demo\proof

\begin{document} 

\begin{abstract}
We apply representation theory to study the homology of
equivariant Dehn-fillings of a given finite, regular cover of a
compact 3--manifold with boundary a torus.  This yields a polynomial
which gives the rank of the part of the homology carried by the solid
tori used for Dehn-filling. The polynomial is a symmetrized form of
the group determinant studied by Frobenius and Dedekind. As a
corollary every such hyperbolic 3--manifold has infinitely many
virtually Haken Dehn-fillings.
\end{abstract}

\maketitle

\section{Introduction} 
In this paper we develop a new connection between the topology of  three-manifolds and the representation theory of finite groups. One consequence for $3$--manifolds is:
\begin{theorem} \label{virthaken}
Let $Y$ be a compact, orientable 3--manifold with one torus boundary component. Assume that the interior of $Y$ admits a complete hyperbolic structure of finite volume. Then infinitely many Dehn-fillings of $Y$ are virtually Haken. 
\end{theorem} 

The case that $Y$ is not fibered also follows from Cooper and Long
\cite{CooperLongVirt}.  The case that $Y$ is fibered but not
semi-fibered is in Cooper and Walsh \cite{CooperWalsh1}. The latter
paper introduced the idea of {\em invariant slope} in a special
case. The investigation of invariant slopes for general finite regular
covers led to the present work.

Unless otherwise stated, in this paper we will use rational coefficients for all homology groups. Throughout $Y$ is a compact, connected, orientable 3--manifold with boundary consisting of a torus $T$ and $\pi\co \tilde{Y}\rightarrow Y$ is a finite regular cover with covering group $G.$ Given a simple closed curve $\gamma$ in $Y$ the {\em virtual rank of $\gamma$} in this cover is the dimension of the subspace of $H_1(\tilde{Y})$ generated by the union of the pre-images of $\gamma.$ This rank only depends on the conjugacy class in $\pi_1(Y)$ determined by  $\gamma.$  
A conjugacy class is a {\em virtual homology class} if there is a finite cover such that this rank is not zero. It has been conjectured that if $Y$ is a hyperbolic $3$--manifold then every non-trivial conjugacy class is a virtual homology class. 
 
  Let $K(\tilde{Y})$ denote the kernel of the map $\incl_*\co H_1(\partial\tilde{Y})\rightarrow H_1(\tilde{Y}).$ For each slope $\alpha$ on $\partial Y$ let $V(\alpha,\tilde{Y})$ denote the subspace of $H_1(\partial\tilde{Y})$ spanned by the pre-images of $\alpha.$  The {\em  filling rank}  of  $\alpha$ for this covering is defined as
$$\mathrm{fillrank}(\alpha,{\tilde{Y}}) =  \dim \left(K(\tilde{Y}) \cap V(\alpha,\tilde{Y})\right).$$
Let $\tilde{Y}(\alpha)$ be the closed 3--manifold obtained by equivariant Dehn-filling of $\tilde{Y}$ along slopes which cover $\alpha.$  Then the filling rank equals the dimension of the subspace, $P,$  of $H_1(\tilde{Y}(\alpha))$ carried by the union of the solid tori used for the Dehn-filling. If the filling rank is greater than zero we say that $\alpha$ is a {\em virtual homology slope} for the covering. The quotient  $\tilde{Y}(\alpha)/G$ is in general an orbifold obtained by an orbifold Dehn-filling of $Y$ and $P$ is the homology carried by the pre-images of the core curve of the (orbifold) solid torus attached to $Y.$ A priori in general one might not expect any interesting virtual homology slopes, and so the following theorem, 
is perhaps surprising:

\begin{theorem} \label{vhs}
Let $Y$ be as in \fullref{virthaken}.  Then
for all $n > 0$, there is a regular cover $\tilde Y \rightarrow Y$ and a slope $\alpha$ on $\partial Y$ so that $\alpha$ is a virtual homology slope for $\tilde Y$ of filling rank bigger than $n.$\end{theorem} 

Suppose $H$ is any subgroup of a finite group $G.$ The collection of sets $S = \{\ HgH\cup Hg^{-1}H\ \}_{g\in G}$ form a partition of $G.$ Define $\sigma\co G\rightarrow S$ by  $\sigma(g) = HgH\cup Hg^{-1}H.$ Introduce commuting variables $Y_s$ for each element $s\in S.$ Enumerate the left cosets of $H$ as $g_1H,\cdots,g_kH$ and define the {\em symmetrized group-coset} matrix $M^{\sym}(G,H)$ to be the $k\times k$ matrix with $(i,j)$ entry   $Y_{\sigma(g)}$  where $g=g_i^{-1}g_j.$ This is related to the group determinant. Here is our main theorem:

\begin{theorem}\label{maintheorem} Suppose that $G$ is a finite group and $Y$ is a compact, connected, orientable $3$--manifold with boundary a torus $T.$ Suppose that $\tilde{Y}\rightarrow Y$ is a regular cover with group of covering transformations $G.$ Let $T_1$ be a component of $\partial\tilde{Y}$ and  $H=\stab(T_1).$ Then there are  $\alpha,\beta$ which generate $H_1(T;{\Bbb Z})$ and are not virtual homology slopes of this cover. Any such basis determines a collection of rationals, $t\co S\rightarrow {\Bbb Q},$ one for each element of $S$ with the following property.

Let $ B$ be the rational matrix obtained from $M^{\sym}(G,H)$ by
setting $Y_s=t(s)$ for each $s\in S.$ Given $m,n\in{\Bbb Z}$ not both
zero then $ \mathrm{fillrank}(m\alpha+n\beta,\tilde{Y})$ equals the
dimension of the $(n/m)$--eigenspace of $ B$.\end{theorem}

The rationals $t(s)$ can be read off from the kernel of the map $\incl_*\co H_1(\partial \tilde{Y})\rightarrow H_1(\tilde{Y}).$ However the power of this theorem lies in the fact that for some pairs $(G,H)$ there are eigenvalues of $M^{\sym}(G,H)$ which are rational linear combinations of the variables. Then every specialization will have certain rational eigenvalues with multiplicity at least as large as for the unspecialized matrix. In \fullref{groupdeterminants} we use representation theory to obtain information about the eigenvalues of $M^{\sym}(G,H).$

As a simple example if $G={\Bbb Z}_3$ and $H=1$ then $\det(M^{\sym}(G,H))=(a-b)^2(a+2b).$ It follows that there are eigenvalues $\lambda=a+2b$ of multiplicity $1$ and $\lambda=a-b$ of multiplicity $2.$ From this we obtain a result which can be easily proved more directly, once it is known.
\begin{corollary}\label{Z3cover} Suppose $Y$ as above and $\tilde{Y}\rightarrow Y$ is a 3--fold cyclic cover and suppose that $T$ has three distinct lifts to this cover. Then there is a virtual homology slope of filling rank $2$  or $3$ 
\end{corollary}

The two possibilities depend on whether or not under the specialization $b=0.$ The corresponding statement for a 5--fold cyclic cover is false.  A more general result, and the proofs of Theorems \ref{virthaken}, \ref{vhs}, and \ref{maintheorem} are given in \fullref{applications}.

\begin{corollary}\label{PSL2cover} Suppose $p$ is a prime and $Y$ as above and $\tilde{Y}\rightarrow Y$ is a regular cover with group of covering transformations $PSL(2,{\Bbb Z}_p).$ Suppose that $T$ is a component of $\partial\tilde{Y}$ and  that $|\stab(T)| = p.$  Then there is a virtual homology slope for this cover of filling rank at least $p.$\end{corollary}

Long and Reid have shown that for a hyperbolic 3--manifold $Y$ there are infinitely many primes for which the hypothesis holds.  The condition $|\stab(T)| = p$ is equivalent to $\stab(T)$ is conjugate to the group of upper-triangular matrices with $1$'s on the diagonal. An explicit example for $n=3$ is given in \fullref{example}. 

A fundamental question which our methods do not address is whether it is possible to have a hyperbolic 3--manifold so that there is a unique slope which is the only virtual homology slope. This does happen for the connect sum of a solid torus and a closed hyperbolic 3--manifold, so hyperbolicity is needed.

In the earlier paper \cite{CooperWalsh1} we studied a certain covering corresponding to $G=H_1(Y,{\Bbb Z}_2).$ The decomposition of $H_1(\tilde{Y})$ into $G$--eigenspaces yields virtual homology slopes of filling rank 1. 
For a more general group $G$ one is led to study the decomposition of $H_1(\tilde{Y};{\Bbb Q})$ into ${\Bbb Q}$--irreducible sub-representations of $G$.

 Sections \ref{linearalgebra} and \ref{groupdeterminants} might be of interest to some algebraists,  are devoid of topology, and can be read independently.
We thank Mark Baker, Marc Lackenby and Alan Reid for helpful conversations. 
The first author was partially supported by NSF grant DMS-0405963.  

\section{Linear algebra}\label{linearalgebra}
To find virtual homology slopes one considers the intersection of the two vector subspaces $K(\tilde{Y})$ and $V(\alpha,\tilde{Y})$ as $\alpha$ varies. These are two subspaces of half dimension in $H_1(\partial\tilde{Y}).$ More generally
we consider the intersection of two subspaces of dimension $n$ in a vector space of dimension $2n.$ One subspace is fixed, but the other rotates. Generically the intersection is dimension zero, but has positive dimension for a finite set of rotation angles. We show below that these angles are determined by a certain polynomial.

Let $V$ be a vector space of even dimension $2n$ over some field ${\Bbb F},$ and let $${\cal B} = (\alpha_1,\alpha_2,\cdots,\alpha_n,\beta_1,\beta_2,\cdots,\beta_n)$$ be an ordered basis of $V.$ Given a point $[x_0:y_0] \in {\Bbb P}^1{\Bbb F}$ in the $1$--dimensional projective space over ${\Bbb F}$ we define an $n$--dimensional subspace of $V$ by
$$V([x_0:y_0],{\cal B}) = Span<\ x_0 \alpha_1 + y_0 \beta_1,\cdots,x_0 \alpha_n+y_0 \beta_n\ >.$$
A non-zero vector in this space is called a {\em vector of constant slope with respect to ${\cal B}$} and $[x_0:y_0]$ is called the {\em slope} of the vector.

\begin{proposition} With $V={\Bbb F}^n\oplus{\Bbb F}^n$ and ${\cal B}$
  as above suppose that $U$ is a half-dimensional subspace of $V.$ Let
  $\pi_i\co {\Bbb F}^n\oplus{\Bbb F}^n\rightarrow{\Bbb F}^{n}$ denote
  projection onto the $i$'th factor. Suppose $S\co {\Bbb
    F}^n\rightarrow{\Bbb F}^n\oplus{\Bbb F}^n$ is a linear map with
  image $U.$ Given $(x,y)\in{\Bbb F}^2\setminus(0,0)$ define $T\co {\Bbb
    F}^n\rightarrow{\Bbb F}^n\oplus{\Bbb F}^n$ by $T(\nu)=(x\cdot
  \nu,y\cdot\nu).$ Set $W(x,y)=U\cap V([x:y];{\cal B}).$ Then:

\begin{itemize}
\item[\rm(i)] $W(x,y) \cong \ker(S\circ\pi_1-T\circ\pi_2).$
\item[\rm(ii)] $U$ contains a vector of constant slope $[x_0:y_0]$ iff $p(x_0,y_0)=0$ where $p(x,y)=\det(S\circ\pi_1 - T\circ\pi_2).$
\item[\rm(iii)] If $U$ contains no vector of constant slope $[0:1]$
  then $W(x,y)$ is isomorphic to the $(y/x)$--eigenspace of
  $(\pi_2\circ S)\circ(\pi_1\circ S)^{-1}.$
\end{itemize}
\end{proposition}

\demo Observe that $W(x,y)$ is the subspace of vectors in $U$ of constant slope $[x:y]$ and that $V([x,y];{\cal B}) = \Im(T).$ Hence
\begin{gather*}
W(x,y) = \Im(S) \cap \Im(T).\\
\tag*{\hbox{Observe that}}
\Im(S) \cap \Im(T) = S\circ\pi_1\left(\ker\left(S\circ\pi_1 - T\circ\pi_2\right)\right).\end{gather*}
The restriction of $S\circ\pi_1$ to $\ker\left(S\circ\pi_1 - T\circ\pi_2\right)$ is injective and so $W(x,y) \cong \ker\left(S\circ\pi_1 - T\circ\pi_2\right)$ which gives (i). 
It follows that 
\begin{gather*}
\dim(W(x,y)) = \nullity\left(S\circ\pi_1 - T\circ\pi_2\right).\\
\tag*{\hbox{Now}}
p(x,y) = \det\left(S\circ\pi_1 - T\circ\pi_2\right)
\end{gather*}
is a homogeneous polynomial of degree $n$ in the variables $x$ and $y.$ It follows that $W(x_0,y_0)$ has positive dimension if and only if  $(x_0,y_0)$ is a root of $p(x,y).$ This proves (ii).  

Observe that
$(\nu_1,\nu_2) \in \ker(S\circ\pi_1 - T\circ\pi_2)$ iff  $$((\pi_1\circ S)\nu_1,(\pi_2\circ S)\nu_1) = (x\cdot \nu_2,y\cdot \nu_2).$$ 
The hypothesis that there is no vector of constant slope $[0:1]$ implies that $\pi_1\circ S$  is injective thus $\nu_1 = x(\pi_1\circ S)^{-1}\nu_2.$ The above is thus equivalent to
$$(\pi_2\circ S)\circ(\pi_1\circ S)^{-1}\nu_2 = (y/x)\nu_2.
\eqno{\qed}
$$

We now describe this in terms of matrices. Choose a basis $(\mu_1,\cdots,\mu_n)$ of $U$ and express these vectors as linear combinations of the basis ${\cal B}$ vectors thus:
$$\mu_j = \sum_{i=1}^n\ (a_{ij} \alpha_i +b_{ij} \beta_i)\qquad\qquad a_{ij},b_{ij}\in F.$$ Define two $n\times n$ matrices $A=(a_{ij})$ and $B=(b_{ij}).$ Let $I_n$ be the $n\times n$ identity matrix then the matrix of $S\circ\pi_1 - T\circ\pi_2$ is the $2n\times 2n$ matrix 
$$M =  \left(\begin{array}{cc} A & -x I_n\\ B & -y I_n\end{array}\right).$$ Observe that the first $n$ columns of this matrix span $U$ and the last $n$ columns span $V([x:y],{\cal B}).$
From the above we see that $p(x,y)=\det(M).$ If $(\nu_1,\nu_2)\in{\Bbb F}^n\oplus{\Bbb F}^n$  then
$$M\cdot
\left(\begin{array}{c}  \nu_1 \\ \nu_2\end{array}\right) 
=  \left(\begin{array}{cc} A & -x I_n\\ B & -y I_n\end{array}\right)
\cdot
\left(\begin{array}{c} \nu_1\\ \nu_2\end{array}\right) 
=\left(\begin{array}{c} A\nu_1-x\nu_2\\ B\nu_1-y\nu_2\end{array}\right).
$$
Thus $(\nu_1,\nu_2)$ is in the kernel of $S\circ\pi_1 - T\circ\pi_2$ if and only if
$$ A\nu_1 = x\nu_2 \quad and\quad  B\nu_1 = y\nu_2.$$
The assumption that there is no vector of constant slope $[0:1]$ implies that  $A$ is invertible. Then  the above is equivalent to
$$ B A^{-1}\nu_2  = (y/x)\cdot\nu_2.$$

It is interesting to contemplate this in the following context. Suppose that $Y$ is a compact manifold with boundary consisting of a collection of $n$ tori, $\{ T_i \}_{i\le i\le n.}$  For each torus we may choose a basis $\alpha_i,\beta_i$ of $H_1(T_i)$ and  apply the proposition  with $$U = \Ker\left[ \incl_*\co H_1(\partial Y)\rightarrow H_1(Y)\right].$$ 

In \fullref{boundarymatrix} we apply this to a finite regular cover $\tilde{Y}$ of a $3$--manifold.  We will see that the basis of $V = H_1(\partial\tilde{Y})$ can be chosen so that $A$ is the identity matrix. The matrix $B$ is then the {\em boundary matrix.} The fact that the intersection pairing vanishes on $U$  implies $B$ is symmetric.  In particular $B$ has real eigenvalues and the algebraic multiplicity of an eigenvalue equals the dimension of the corresponding subspace. We will see that the polynomial $p(x,y)$ is related to the {\it group determinant.}

\section{Virtual homology slopes and the boundary matrix}\label{boundarymatrix}

Given a regular cover $\tilde{Y}\rightarrow Y,$ we define a matrix of rational numbers, the {\em boundary matrix,} which encodes $K(\tilde{Y}).$ This matrix has certain symmetry properties. We also define another matrix, {\em the general boundary matrix,} whose entries are variables satisfying the same symmetry properties. The boundary matrix is obtained by replacing the variables in the general boundary matrix by certain rationals. In \fullref{groupdeterminants} we show that the general boundary matrix is closely connected to the group determinant.

We define a {\em slope} on a torus $T$ to be an element of the projective space ${\Bbb P}H_1(T;{\Bbb Q}).$ A slope is uniquely determined by any of the following:  an essential simple closed curve on $T,$ a primitive element of $H_1(T;{\Bbb Z})$ or a non-zero element of $H_1(T;{\Bbb Q}).$  In each case there is an equivalence relation that gives a bijection between equivalence classes and slopes. We will find it convenient to suppress mention of this relation.  This does not lead to any serious ambiguity.

Throughout this section $Y$ is a compact oriented 3--manifold with boundary a torus $T$ and  we fix a regular covering $\pi\co \tilde{Y}\rightarrow Y$ with covering group $G.$ Let $T_1,T_2,\cdots,T_k$ denote the boundary components of $\tilde{Y}.$  The orientation on $Y$ induces ones on $T, \tilde{Y}$ and $\partial\tilde{Y}.$ We will use $\iota(\_\   , \_)$ to denote the (skew-symmetric) intersection pairing on both $H_1(\partial Y;{\Bbb Q})$ and  $H_1(\partial\tilde{Y};{\Bbb Q}).$ 
\begin{lemma}\label{kvhslopes} There are at most $k$ distinct slopes on $Y$ that are virtual homology slopes for the covering $\tilde{Y}.$
\end{lemma} 
\demo  We assign a vector, $v=v(\alpha)\in{\Bbb Q}^k$ to each virtual homology slope $\alpha\in H_1(\partial Y;{\Bbb Q})$ as follows. Define $\alpha_i\in H_1(T_i;{\Bbb Q})$ by $\pi_*\alpha_i=\alpha.$ There is  a non-zero $v=(v_1,v_2,\cdots,v_k)\in{\Bbb Q}^k$ so that $\sum_i v_i\alpha_i$ is in $K(\tilde{Y}).$ The lemma follows from the fact that the vectors assigned to distinct virtual homology slopes are orthogonal. To see orthogonality,  suppose that $\beta$ is another virtual homology slope and $\beta_i\in H_1(T_i;{\Bbb Q})$ satisfies $\pi_*\beta_i=\beta.$ Since the cover is regular it follows that $\iota(\alpha_j,\beta_j)$ is independent of $j.$ Also $\iota(\alpha_i,\beta_j)=0$ if $i\ne j$ since $\alpha_i,\beta_j$ are on different tori. Let $u$ be the vector assigned to the virtual homology slope $\beta.$ The intersection pairing vanishes on $K(\tilde{Y})$ thus 
$$0 = \iota(\sum_i v_i\alpha_i , \sum_j u_j\beta_j) = \iota(\alpha_1,\beta_1)\cdot \sum_iv_iu_i.$$
 Since $\alpha_1$ and $\beta_1$ are distinct slopes on a torus they
 have non-zero intersection number, thus $\sum_i v_iu_i=0$. This
 proves orthogonality. \endproof
 
An ordered basis $(\alpha,\beta)$ of $H_1(T;{\Bbb Z})$ is called {\em
generic} for $\tilde{Y}$ if neither $\alpha$ nor $\beta$ is a virtual
homology slope for this cover. We fix such a generic basis with
$\iota(\alpha,\beta) = +1.$ This gives an identification ${\Bbb
P}H_1(T;{\Bbb Q})\equiv {\Bbb Q}\cup\infty.$ The slope
$m\alpha+n\beta$ corresponds to $n/m\in{\Bbb Q}\cup\infty.$ The fact
that $\alpha$ and $\beta$ are not virtual homology slopes means that
every virtual homology slope lies in ${\Bbb Q}\setminus 0$ with
respect to this basis.  Let $(\alpha_i,\beta_i)$ be the ordered basis
of $H_1(T_i; {\Bbb Q})$ which projects by $\pi_*$ to $(\alpha,
\beta)$.  Observe that $\iota(\alpha_i,\beta_j) = d^{-1}
\delta_{ij}$ where $d=\mathrm{degree}(\pi|:T_1\rightarrow T) = |G|/k.$
   
 \begin{lemma} There is a unique $\gamma_1\in K(\tilde{Y})$ and $b_1,\cdots,b_k\in{\Bbb Q}$ such that
 $$\gamma_1 =\alpha_1 + \sum_{j=1}^k b_j\beta_j.$$
  \end{lemma}
  \demo We have $H_1(\partial\tilde{Y}) = V(\alpha) \oplus V(\beta)$
  where $V(\alpha)$ is the subspace spanned by
  $\alpha_1,\cdots,\alpha_k$ and $V(\beta)$ is similarly
  defined. Consider the projection  onto the first factor
  $p\co H_1(\partial\tilde{Y}) \rightarrow V(\alpha).$ Then $p$
  restricted to $K(\tilde{Y})$ is injective, otherwise $\beta$ is a
  virtual homology slope. Since $K(\tilde{Y})$ and $V(\alpha)$ have
  the same dimension it follows that this restriction is an
  isomorphism. Hence there is a unique $\gamma_1$ which projects onto
  $\alpha_1.$
\endproof

A compact, connected, oriented surface $S$ properly embedded in $\tilde{Y}$ for which $[\partial S]$ is a non-zero multiple of $\gamma_1$ is called a {\em fundamental surface.} It depends on the choice of generic basis of $H_1(\partial Y)$ and on the labeling of the components of $\partial\tilde{Y}.$ Such a surface always exists.

By applying covering transformations, it follows that for each torus $T_i$ there is a unique  $\gamma_i\in K(\tilde{Y})$ such that
$$\gamma_i = \alpha_i + \sum_{j=1}^k b_{ij}\beta_j.$$ 
 We call the $k\times k$ matrix $B = (b_{ij})$ the {\em boundary matrix}  for $\tilde{Y}.$ The boundary matrix depends on the choice of generic basis and on the ordering of the components of $\partial\tilde{Y}.$
It follows that 
$$b_{ij}\ =\ d\cdot \iota( \alpha_j,\gamma_i ).$$


\begin{lemma}\label{Bissymmetric}  The boundary matrix $B$ is invertible and symmetric.
\end{lemma}  
\demo If the boundary matrix is singular then some linear combination of the rows of $B$ are zero. The corresponding linear combination of the $\gamma_i$'s  is a non-zero element $\gamma\in K(\tilde{Y})$ which is a linear combination of $\alpha_1,\cdots,\alpha_k$.  This contradicts that $\alpha$ is not a virtual homology slope for $\tilde{Y}$. To prove symmetry we use the fact that 
the intersection pairing vanishes on $K(\tilde{Y})$ and recalling that
$\iota(\alpha_i, \beta_j) = d^{-1} \delta_{ij}$ gives
$$\eqalignbot{
0 &= \iota( \gamma_i,\gamma_j )\cr
&=  \iota(\alpha_i + \sum_k b_{ik}\beta_k \qua ,\qua   \alpha_j + \sum_l b_{jl}\beta_l )\cr
&= d^{-1} (-b_{ji} + b_{ij}).}
\eqno{\qed}$$

The symmetry may be seen another way. Assume the slope $\beta\subset \partial Y$ lifts to $\tilde{Y}$.  Then the Dehn-filling, $Y(\beta),$ of $Y$ along $\beta$ is covered by a Dehn-filling, $\tilde{Y}(\beta),$ of $\tilde{Y}$ along lifts of $\beta.$ For simplicity, suppose that $\tilde{Y}(\beta)$  is a homology 3--sphere. Then there is a surface $S_i$ in $\tilde{Y}(\beta)$ with boundary  which is the core curve, $\gamma_i,$ of the $i$'th solid torus in $\tilde{Y}(\beta).$ The surface $S_i$ can be chosen so it meets $\gamma_j,$ for every $j\ne i,$ transversally and minimally. Thus $S_i\cap\tilde{Y}$ has boundary some longitude $ \alpha_i\beta_i^n$ of $T_i$ together with meridians $ \beta_j$ on the other $T_j$ for $j\ne i.$

The {\em linking number} of $\gamma_i$ and $\gamma_j$ in $\tilde{Y}( \beta)$ is 
$$Lk(\gamma_i,\gamma_j) = \#(S_i\cap\gamma_j) = b_{ij}.$$
Now the linking number is {\em symmetric}, ie $Lk(\gamma_i,\gamma_j) = Lk(\gamma_j,\gamma_i)$, which shows again that $B$ is a {\em symmetric} matrix.

Let $H$ be the subgroup of $G$ which stabilizes $T_1.$ For $g\in G$  the subset of $G$ consisting of elements which map $T_1$ to $g\cdot T_1$ is $g\cdot H.$ Thus components of $\partial\tilde{Y}$ are in one to one correspondence with left cosets of $H$ in $G.$ For each $1\le i\le k$ choose $g_i\in G$ so that $g_i\cdot T_1=T_i.$ Then  $\lbrace g_1,\cdots,g_k \rbrace$ is a complete set of left coset representatives of the subgroup $H$ of $G.$
  
  \begin{corollary}  The $G$--orbit of $\gamma_1$ spans $K(\tilde{Y}).$ In particular $K(\tilde{Y})$ is a cyclic representation of $G$ which is isomorphic to the action of $G$ by permutations on left cosets of $H$ in $G.$
  \end{corollary}
  \demo The elements $\gamma_1,\cdots,\gamma_k$ are linearly independent in $K(\tilde{Y}).$ Since $\dim(K(\tilde{Y}))$ $= k$ they form a basis.   Since $g_iT_1=T_i$ it follows that $g_i\gamma_1=\gamma_i$ so the $G$--orbit of $\gamma_1$ spans $K(\tilde{Y}).$ Since $\gamma_i\subset T_i$ the action of $G$ on $K(\tilde{Y})$ is isomorphic to the action of $G$ by permutations on the components of $\partial\tilde{Y}.$ This, in turn, is isomorphic to the action of $G$ on left cosets of $H.$
  \endproof

\begin{proposition}[Virtual homology slopes are eigenvalues of the boundary matrix]\label{slopesareeigenvalues}\  \\ The slope $t\in{\Bbb Q}$  (which corresponds to $\alpha+t\beta$) is a virtual homology slope of the covering $\tilde{Y}\rightarrow Y$ if and only if $t$ is an eigenvalue of the boundary matrix $B.$ Furthermore the filling rank of this slope equals the dimension of the $t$--eigenspace of $B.$\end{proposition}
\demo First suppose that $\alpha+t\beta$ is a virtual homology
slope. Then there is a non-zero element of $K(\tilde{Y})$ of the form
$$\gamma = \sum_{i=1}^k\ c_i(\alpha_i + t\beta_i).$$
Now $\gamma$ is a linear combination of the basis elements
$\gamma_1,\cdots\gamma_k$ of $K(\tilde{Y}).$ Recalling that
$\gamma_i=\alpha_i+(\beta's)$ we get that   
\begin{gather*}
\gamma\ = \sum_{i=1}^k\ c_i\gamma_i\ =\ \sum_{i=1}^k\
c_i[\alpha_i+\sum_{j=1}^k b_{ij}\beta_j].\\
\tag*{\hbox{The part of $\gamma$ on  torus $T_j$ is}}
c_j\alpha_j \ +\ \sum_{i=1}^k b_{ij}c_i\beta_j.\\
\tag*{\hbox{This has slope $t$ so we get}}
\sum_{i=1}^k b_{ij}c_i\ =\ tc_j.
\end{gather*}
Recalling that $B$ is symmetric this implies
\begin{gather*}
\sum_{i=1}^k b_{ji}c_i\ =\ tc_j.\\
\tag*{\hbox{In other words}}
B\vec{c}\ = t\cdot\vec{c}
\end{gather*}
where $\vec{c}=(c_1,c_2,\cdots,c_k).$ Thus $\vec{c}$ is an eigenvector of $B$ with eignevalue $t.$ This also shows there is an injective linear map from the subspace, $W(t),$ of $K(\tilde{Y})$ consisting of vectors of slope $t$ into the $t$--eigenspace, $E_t,$ of $B.$ Hence $\dim(W_t)\le \dim(E_t).$

For the reverse inequality suppose that $\vec{c}=(c_1,c_2,\cdots,c_k)\in E_t$ so $B\cdot\vec{c} = t\cdot\vec{c}.$ Then $$\gamma = \sum_{i=1}^k\ c_i  \gamma_i $$ is in $K(\tilde{Y}).$
 Using the calculations above we see that the components of $\gamma$ on each torus $T_j$ all have slope $t.$ It follows that there is an injective linear map from $E_t$ into $W(t).$ Hence $W(t)$ and $E_t$ have the same dimension.
 \endproof

The boundary matrix may be regarded as a function on pairs of left $H$--cosets, $B\co G/H\times G/H\rightarrow{\Bbb Q}$ given by $$B(g_iH,g_jH) = b_{ij}.$$

\begin{proposition}\label{boundarymatrixidents}  
The boundary matrix  satisfies the following properties:
\begin{enumerate}
\item Invariance under left translation: $B(xH,yH) = B(gxH,gyH)$
\item Symmetry: $B(xH,yH)=B(yH,xH)$
\item Inversion: $B(H,g^{-1}H) = B(H, gH)$
\end{enumerate}
\end{proposition}
\demo
Recall that 
\begin{gather*}
b_{ij}\ =\  d\cdot\iota( \alpha_i,\gamma_j ).\\
\tag*{\hbox{Thus}}
B(xH,yH)\ =\ d\cdot\iota( x\alpha_1,y\gamma_1 ).
\end{gather*} 
Intersection numbers are preserved by covering transformations, thus left invariance follows from $\iota( x\alpha_1,y\gamma_1 ) = \iota( gx\alpha_1,gy\gamma_1 ).$
 Symmetry follows from \eqref{Bissymmetric}. Finally using the first two properties we get $B(H,g^{-1}H) = B(gH,H) = B(H,gH)$ which gives (3).\endproof

It follows that $B(xH,yH)=B(H,x^{-1}yH)$ so the boundary matrix is determined by its values in the first row, ie by the function $L\co G/H\rightarrow{\Bbb Q}$ given by $L(xH)=B(H,xH).$ We describe this by saying $B$ is {\em induced} by $L.$ In other words the entries in the boundary matrix are not arbitrary but satisfy certain relations of the form $(i,j)$--entry equals $(k,l)$--entry for values of $i,j,k,l$ that only depend on $G$ and $H.$ Thus the boundary matrix is determined by $\gamma_1$ together with how $G$ permutes the entries in the first row of $B$ to give the other rows of $B.$ Geometrically this corresponds to the statements that covering translates of the fundamental surface generate $K(\tilde{Y})$ and the boundary of the fundamental surface is given by the first row of the boundary matrix. 

\begin{lemma} A function $L\co G/H\rightarrow{\Bbb Q}$ induces a function $B\co G/H\times G/H\rightarrow{\Bbb Q}$ which satisfies properties (1) to (3) iff $L$ factors through a map $\overline{L}\co S\rightarrow{\Bbb Q}$ where $$S = \{\ HgH\cup Hg^{-1}H\ \}_{g\in G}$$ is the partition of $G$ defined in the introduction.\end{lemma}
\demo First suppose we are given $\overline{L}\co S\rightarrow{\Bbb Q}$ then we obtain a well defined function $B\co G/H\times G/H\rightarrow{\Bbb Q}$ given by $B(xH,yH)=\overline{L}(Hx^{-1}yH\cup Hy^{-1}xH).$ It is immediate that $B$ satisfies (1), (2) and  (3).

For the converse, given $B\co G/H\times G/H\rightarrow{\Bbb Q}$ which
satisfies properties (1) to (3) define $L\co G/H\rightarrow{\Bbb Q}$ as
above. Suppose that $HxH\cup Hx^{-1}H = HyH\cup Hy^{-1}H.$ Then either
$y\in HxH$ or $y\in Hx^{-1}H.$ In either  case $yH=hx^{\epsilon}H$ for
some $h\in H$ and $\epsilon=\pm1.$ Using properties (1) to (3) for $B$
we get 
\begin{align*}
L(yH)=B(H,yH)=B(H,hx^{\epsilon}H)&=B(h^{-1}H,x^{\epsilon}H)\\
&=B(H,x^{\epsilon}H)=B(H,xH)=L(xH).\end{align*}
It follows that $L$ factors through $\overline{L}.$\endproof

{\bf Remark}\qua A double coset $HgH$ may be viewed as the union of left cosets $xH$ where $x = hgh^{-1}$ varies over the $H$--conjugates of the element $x\in G.$ Thus the pair of double cosets $HgH \cup Hg^{-1}H$ can be viewed as the union of left $H$ cosets obtained by taking $H$--conjugates of an unoriented loop (ie pair of elements $g$ and $g^{-1}$).

The  $(i,j)$ entry of the boundary matrix is
$$b_{ij} = B(g_iH,g_jH)  =  B(H,g_i^{-1}g_jH) = L(g_i^{-1}g_jH).$$ This is determined by  which element of $S$ contains $g_i^{-1}g_j.$  Let $\{\ Y_s\ |\ s\in S\ \}$ be a set of commuting variables, one for each element of $S.$ Recall that $g_1,\cdots,g_k$ is a complete set  of left coset representatives of the subgroup $H$ of $G$ and $\sigma\co G\rightarrow S$ is given by  $\sigma(g) = HgH\cup Hg^{-1}H.$  The {\em general boundary matrix} ${\cal B}$  is the $k\times k$ matrix whose $(i,j)$ entry is $Y_{\sigma(g_i^{-1}g_j)}.$ Observe that this matrix only depends on the pair $(G,H)$ together with an ordering of the left cosets of $H$ (coming from the ordering of the boundary components of $\partial\tilde{Y}.$) {\em In fact we see from the definition given in the introduction that ${\cal B}$ equals $M^{\sym}(G,H).$}  In particular, when the boundary torus lifts, the generalized boundary matrix is the symmetrized group matrix defined in \fullref{groupdeterminants}. The boundary matrix for a particular 3--manifold is obtained by replacing the variables in the general boundary matrix by certain rationals, ie by a specialization $t\co\{\ Y_s\ |\ s\in S\ \}\rightarrow {\Bbb Q}$ of the variables. We record this for later use:

\begin{theorem}\label{genbdry} Suppose $Y$ is a compact oriented $3$--manifold with boundary a torus $T$ and that $\tilde{Y}\rightarrow Y$ is a regular cover with $G$ the group of covering transformations. Let $T_1,\cdots,T_k$ be the boundary components of $\tilde{Y}$ and let $H$ be the stabilizer of $T_1.$ Choose a generic basis of $H_1(T)$ and let $B$ be the boundary matrix defined with this data. Then
there is a specialization $t\co \{\ Y_s\ |\ s\in S\ \}\rightarrow {\Bbb Q}$ of the variables so that $B = t(M^{\sym}(G,H)).$\end{theorem}

 In \fullref{groupdeterminants} we show how the symmetrized group-coset matrix is related to group determinants and representation theory. Although we will not make use of this fact,  every specialization of the general boundary matrix is a boundary matrix for some covering of some 3--manifold:

\begin{proposition} Let $G$ be a finite group and $H < G$ such that $\mathbb{Z} \oplus \mathbb{Z}$ surjects $H$. With the above notation, given $t\co \{\ Y_s\ |\ s\in S\ \}\rightarrow {\Bbb Q}$ there is a $3$--manifold $P$ and a regular cover $\tilde{P}$ of $P$ with covering group $G$ such that the boundary matrix is $t(M^{\sym}(G,H)).$\end{proposition}
\demo We will only sketch the proof in the case that $H=1.$ 
We will  construct a manifold $\tilde{P}$ with a free $G$--action and set $P=\tilde{P}/G.$ The fact that $H=1$ means that $|\partial \tilde{P}| = |G| = k$ and $G$ freely permutes the components of $\partial\tilde{P}.$ Let $m>0$ be an integer so that $C_{ij} = m\cdot t(\sigma(g_i^{-1}g_j))$ is an integral matrix. Set $a=\gcd(m,|C_{11}|)$ and $b=\sum_{i=2}^{k} |C_{1i}|.$ 

The manifold $\tilde{P}$ we seek has boundary consisting of tori $T_1,\cdots T_k$ where $T_i=g_i\cdot T_1.$ It  contains a fundamental surface $F$ with $a+b$ boundary components of which $b$ boundary components are on  $T_2,\cdots T_k$ and project to loops parallel to $\beta.$ In fact $T_i\cap F = C_{1i}\beta_i\in H_1(T_i)$  for $i\ge 2.$ The remaining $a$ boundary components are parallel loops on $T_1$ whose sum is the homology class $m\alpha_1+C_{11}\beta_1\in H_1(T_1).$  The orbit of the fundamental surface under $G$ is  a collection of surfaces corresponding to the rows of $C.$ The number of boundary components of $g_i\cdot F$ on the torus $g_j\cdot T_1$ is $|C_{ij}|$ for $i\ne j,$ and $a$ for $i=j,$ and the sign of $C_{ij}$ determines the direction a boundary component of $F$ winds round $T.$ The existence of these surfaces in $\tilde{P}$ ensures that the boundary matrix is a scalar multiple of $C.$

These surfaces may be assumed to be transverse to one another. They intersect along circles and arcs which run between intersection points on the boundary of opposite sign. We will assume there are no circles and that the number of such arcs is minimal. The latter is just the requirement that the various loops which are the boundary components of the fundamental surfaces have minimal intersection with each other. Given $g_i\ne g_j\in G$ the boundaries of the surfaces $F_i = g_i\cdot F$ and $F_j =g_j\cdot F$  only intersect on $T_i$ and $T_j.$ This is because on the remaining tori they both have boundary components parallel to $\beta.$ Minimality of intersection then implies they do not intersect on any other torus. Furthermore minimality also implies all the intersections on one torus have the same sign. Thus these two surfaces intersect along arcs each of which has one endpoint on each of these two tori. The number of such arcs is given by the absolute value of the algebraic intersection of $F_i\cap T_i$ with $F_j\cap T_i$ and can be computed from $C.$ 

We decompose  $\tilde{P} = N\cup (\tilde{P}\setminus int(N))$ where $N$ is a regular neighborhood of $\partial\tilde{P} \cup \bigcup F_i.$ To construct $\tilde{P}$ we reverse the above. We construct $N$ by gluing copies of $F$ onto  $\partial\tilde{P}$ and then thickening. The construction must be done $G$--invariantly.  For each $g\in G$ there is one copy, $F_g,$ of $F$ and a torus $T_g.$ The group acts on this set by permuting the labels. It is easy to use $C$ to $G$--invariantly glue each $F_g$ to the tori as described above. We need a $G$--invariant collection of disjoint arcs in $\cup F_i$ to enable the gluing together of the copies of $F.$ On $F_1$ choose a set of disjoint arcs, one going from each  component of $F_1\cap T_1$ to each component of $F_1\cap T_i$ for $i\ge 2.$ We can certainly do this if $F$ has large enough genus. Now take the required number of parallel copies of these arcs. Using $G$ we can move copies of these arcs to each $F_i.$

This allows us to construct a $2$--complex by gluing the copies of $F$ and the tori in the way required by $G.$ This complex can be thickened to be a $3$--manifold $N$ on which $G$ acts freely. The boundary of $N$ is $\cup T_g$ together with some other boundary components $S$ coming from parts of the neighborhood of copies of $F.$ Take two copies of $N$ and glue them along $S.$ The result is a manifold with a free $G$--action and boundary two copies of $\cup T_g.$ Equivariantly Dehn-fill along one copy of $\cup T_g.$ This produces a compact manifold $\tilde{P}$ with a free $G$--action, boundary $\cup T_g,$ and the required boundary matrix.\endproof

\section{Group determinants}\label{groupdeterminants}

The theory of group determinants was initiated by Frobenius and
Dedekind. An interesting history of this and the origins of
representation theory is contained in Lam \cite{Lam}. A contemporary
survey is Johnson \cite{Jo}. We introduce a version of the group determinant
that is {\em symmetric relative to a subgroup.} We also show how the
symmetrized group-coset matrix defined in the introduction is related
to this.

In what follows, for simplicity we will work with complex vector spaces. Suppose $G$ is a group, $V$ is a finite dimensional vector space  and $\rho\co G\rightarrow \Aut(V)$ is a  representation.  Consider the polynomial ring $\Lambda = {\Bbb C}[\{\ X_g\ :\ g\in G\ \}]$ with one indeterminate, $X_g,$ for each element $g$ of the group $G.$ The {\em representation matrix} of $\rho$ is the matrix with entries in this ring $$M(\rho)\ =\  \sum_{g\in G} \rho(g)\cdot X_g \ \in\ M_n(\Lambda).$$ One may view this as the image of the ``general element'', $\sum_g X_g\cdot g \in \Lambda[G],$  of the group $G$ under the induced ring homomorphism $\rho_*\co \Lambda[G]\rightarrow \End(V\otimes \Lambda).$
The {\em representation determinant}  of $\rho$ is   $$\det(\rho) = \det\left(M(\rho)\right)$$ which is a homogeneous polynomial in  $\Lambda$ of degree $\dim(V).$
 In particular if $\rho$ is a $1$--dimensional representation, then $\det(\rho)$ is a polynomial of degree $1.$
The {\em right regular representation,} $R_G\co G\rightarrow{\Bbb C}^{|G|}$ is the representation obtained from the action of $G$ by multiplication on the group ring ${\Bbb C}[G]$ induced by the right regular action of $G$ on itself, $h\co  g \mapsto gh^{-1}$.  A basic result in the representation theory of finite groups is:
\begin{theorem}\label{rightreg} The right regular representation is conjugate to
$$\bigoplus_{\rho}\left(\bigoplus_{i=1}^{\dim(\rho)} \rho\right)$$ where the first sum is over all conjugacy classes of irreducible representations.\end{theorem} The {\em group matrix}  of $G$ is $M(G) = M(R_G).$  The group matrix was originally defined using the left-regular representation, but the theory is exactly the same. Let $\{\ g_i\  :\ 1\le i\le n\ \}$ be a list of the elements of $G.$  Another description of this matrix is that it is the $n\times n$ matrix (where $n=|G|$) whose $(i,j)$ entry is $X_{g_i^{-1}g_j}$.  This follows since the right-regular action of $g_i^{-1}g_j$ takes $g_j$ to $g_i$.  The {\em group determinant} of $G$ is the representation determinant of the right regular representation, ie $$\det(G) = \det\left(M(G)\right).$$
Since the right regular representation involves integral matrices, this polynomial has integer coefficients. The following is immediate:

\begin{proposition}\label{decompose} $M(\rho_1\oplus\rho_2)=M(\rho_1)\oplus M(\rho_2)$ and $\det(\rho_1\oplus\rho_2) = \det(\rho_1)\cdot \det(\rho_2).$\end{proposition}

This is useful for computing $\det(G),$ since instead of having to work with the determinant of an $n\times n$ matrix, where $n=|G|,$ coming from the right regular representation, one can instead compute the contributions to the group determinant from each irreducible representation separately. The size of the matrices involved is then the dimension of the irreducible representation. Since the regular representation of an abelian group is the sum of the one-dimensional representations one obtains:

\begin{corollary}\label{abelian}
The group determinant of a finite abelian group is
$$\det(G) = \prod_{\rho\co  G \rightarrow \mathbb{C}^*}(\sum_{g \in G} \rho(g)\cdot X_g)$$
\end{corollary}

\noindent Given a subgroup $H$ of $G$ we form the following quotients
of $\Lambda$
\def\struttt{\vrule width 0pt depth 7pt}
$$\begin{array}{lll}
\pi^{\sym}: & \Lambda\rightarrow\Lambda^{\sym}\ & =\ \Lambda / < X_g-X_{g^{-1}}\ |\ g\in G\ >\struttt\\
\pi_H^{\phantom{\sym}}: & \Lambda\rightarrow\Lambda_H\ & =\ \Lambda / < X_g-X_{gh}\ |\ g\in G,\ h\in H\  >\struttt\\
\pi^{\sym}_H :& \Lambda\rightarrow\Lambda^{\sym}_H\ & =\ \Lambda / < X_g-X_{g^{-1}}\ ,\ X_g-X_{gh}\ |\ g\in G,\ h\in H\  >
\end{array}
$$
Each of these quotients is a polynomial ring obtained by formally identifying some of the variables in $\Lambda.$ Generators of $\Lambda^{\sym}$ correspond to subsets $\{g,g^{-1}\}\subset G$ and generators for $\Lambda_H$ correspond to left cosets of $H.$ Generators of $\Lambda^{\sym}_H$ correspond to the subsets $HgH\cup Hg^{-1}H\subset G.$ Observe that the ideal used for the third quotient is the ideal generated by the previous two ideals. 

The {\em symmetrized group determinant} is $${\det}^{\sym}(G)\ =\
\pi^{\sym}(\det(G)).$$ Sjogren introduced this in connection with
groups acting on graphs (Sjogren \cite{Sj}).  Informally one just equates
$X_g\equiv X_{g^{-1}}$ for all $g\in G$ in the group determinant.  The
{\em symmetrized group matrix} of $G$ is $$M^{\sym}(G) =
\pi^{\sym}(M(G)).$$ From the description of the $(i,j)$ entry of the
group matrix as $X_{g_i^{-1}g_j}$ it follows that this is a symmetric
matrix. Hence the symmetrized group determinant is the determinant of
this symmetric matrix. Observe that if every element of $G$ is its own
inverse then the symmetrized group determinant equals the group
determinant.

Here are some examples. In the sequel we will be interested in how many linear factors there are over ${\Bbb Q}$, and what the multiplicity of these factors is, in the symmetrized group determinant.

\begin{enumerate}
\item
For $G = {\Bbb Z}_3$ the group determinant is
$$\det({\Bbb Z}_3)\ =\ (a + b + c)(a^2 - a b + b^2 - a c - b c + c^2)$$ where the correspondence between group elements and variables is $(0,1,2)\leftrightarrow(a,b,c).$ If we  symmetrize this by setting $b=c$ we obtain
$${\det}^{\sym}({\Bbb Z}_3)\ =\ (a + 2 b)(a - b)^2 .$$

\item  ${\det}^{\sym}({\Bbb Z}_4)\ =\ (a+2b+c)(a-2b+c)(a-c)^2 .$

\item  ${\det}^{\sym}( {\Bbb Z}_5)\ =\  (a + 2 b + 2 c)(a^2 - a b - b^2 -
       a c + 3 b c - c^2)^2 .$
       
\item  ${\det}^{\sym}( {\Bbb Z}_6)\ =\ (a + 2b + 2c + d)(a + b - c - d)^2 (a - 2 b + 2 c - d)(a - b - c + d)^2 .$

\item   ${\det}^{\sym}( {\Bbb Z}_8)\ =\ (a + 2 b + 2 c + 2 d + e)(a - 2\ c + e)^2(a - 2 b + 2 c - 2 d + e)(a^2 - 2 b^2 + 4 b d - 
      2 d^2 - 2 a e + e^2)^2 .$
      
\item The dihedral group with six elements:\\
${\det}^{\sym}( D_6)\ =\ (a + 2 b + d + e + f)(a +  2 b - d - e - f)(a^2 - 2 a b + b^2 - d^2 + d e - e^2 +  d f + e f - f^2)^2 .$

\item The quaternionic group with 8 elements:\\
${\det}^{\sym}(Q_8)\ =\ (a + 2 b + 2 c + 2 d + e)(a - e)^4(a + 2 b - 2 c - 2 d + e)(a - 2 b + 2 c - 2 d + e)(a - 2 b - 2 c + 2 d + e) .$

\item The Klein four group:\\
${\det}^{\sym}({\Bbb Z}_2\oplus{\Bbb Z}_2)\ =\ (a + b + c + d)(a + b - c - d) (a - b + c - d) (a - b - c + d) .$

\item ${\det}^{\sym}({\Bbb Z}_3\oplus{\Bbb Z}_3)\ =\ (a + 2b +2d + 2e +2f)(-a +b +d+e -2f)^2(a -b + 2d -e -f)^2 (-a +b +d -2e +f)^2 (-a - 2b +d +e +f)^2 .$
\end{enumerate}

It is a result of Frobenius that the irreducible factors over ${\Bbb C}$ of the group determinant are the representation determinants of the irreducible representations appearing in the right regular representation. We do not know the general factorization of the symmetrized group determinant over ${\Bbb C}$ or ${\Bbb Q}.$

 Recall the definition of the symmetrized group-coset matrix given in the introduction. Given a subgroup $H$ of $G$ enumerate the left cosets as $g_1H,g_2H,\cdots,g_kH$ with $g_1H=H.$  Consider the matrix whose $(i,j)$--entry is $\pi^{\sym}_H(X_{g_i^{-1}g_j}) .$ This entry is $X_{HgH\cup Hg^{-1}H}$ where $g=g_i^{-1}g_j.$ If we identify this variable with $Y_{\sigma(g)}$ then this matrix is just $M^{\sym}(G,H).$ 

If $H$ is the trivial subgroup then this is just the symmetrized group matrix. In general the symmetrized group-coset matrix is a symmetric non-singular matrix. Indeed specializing all the variables equal to zero except $Y_H$ one obtains $Y_H$ times the identity matrix.

\begin{proposition} \label{groupcoset} Choose a complete set of left coset representatives $g_1,\cdots g_k$ of $H$ in $G.$ Let $B=M^{\sym}(G,H)$ be the symmetrized group-coset matrix for $(G,H)$ with respect to this ordered basis. Enumerate the elements of $H$ as $h_1,\cdots,h_l$ then enumerate the elements of $G$ so that for $1\le p\le k$ and $1\le q\le l$ element $p+(q-1)k$ is $g_ph_q.$  With this as an ordered basis of ${\Bbb C}[G]$,     $$\pi^{\sym}_H(M(G))\   = \ \left(\begin{array}{cccc}
B & B & \cdots & B\\
B & B & \cdots & B\\
\cdots & \cdots &\cdots &\cdots\\
B & B & \cdots & B\\
\end{array}\right).$$
Furthermore $\pi^{\sym}_H(M(G))$ is conjugate to $ l M^{\sym}(G,H)\oplus 0.$ 
It follows that for $\lambda\ne0$ the dimension of the $\lambda$--eigenspace of $M^{\sym}(G,H)$  equals the dimension of the $(l\lambda)$--eigenspace of $\pi^{\sym}_H(M(G)).$\end{proposition}
\demo The $(i,j) = (p+(q-1)k,r+(s-1)k)$ entry of $\pi^{\sym}_H(M(G))$ is $\pi^{\sym}_H(X_g)$ where $$g = (g_p h_q)^{-1}(g_r h_s)= h_q^{-1}(g_p^{-1}g_r)h_s.$$ The $(p,r)$ entry of $M^{\sym}(G,H)$ is $\pi^{\sym}_H(g_p^{-1}g_r) .$ This equals $\pi^{\sym}_H(X_g)$ for $g=g_p^{-1}g_r.$

For the second conclusion, there is an invertible matrix $Q$ with entries in ${\Bbb Q}$ so that
$$ Q \ \left(\begin{array}{cccc}
B & B & \cdots & B\\
B & B & \cdots & B\\
\cdots & \cdots &\cdots &\cdots\\
B & B & \cdots & B\\
\end{array}\right) Q^{-1} = \left(\begin{array}{cccc}
l B & 0 & \cdots & 0\\
0 & 0 & \cdots & 0\\
\cdots & \cdots &\cdots &\cdots\\
0 & 0 & \cdots & 0\\
\end{array}\right)\eqno{\lower25pt\hbox{\qed}}$$

The following shows how the eigenvalues of the symmetrized group-coset matrix can be determined from the roots of the group determinant.

\begin{proposition}\label{grpcosetgrpdet} The characteristic polynomial of $M^{\sym}(G,H)$ can be obtained from the group determinant, $\det (G),$ of $G$ by a change of variables. For $g\ne1$ replace $X_g$ by $Y_{\sigma(g)}.$ Replace $X_1$ by $Y_H - t.$ The resulting polynomial is
$t^m p(t)$ where $p(l t)$ is the characteristic polynomial of $M^{\sym}(G,H)$ and $m = |G| - |H|.$
\end{proposition}
\demo Let $n=|G|$ and  $I_n$ be the $n\times n$ identity matrix. Then  $M(G) = X_1I_n + P$ where $P$ is a matrix that does not involve $X_1.$ Thus the  characteristic polynomial of $M(G)$ is obtained from $\det(G)$ by replacing $X_{1}$ with $X_{1} - t$.  It follows that the characteristic polynomial of $\pi^{\sym}_H(M(G))$ is obtained from this by now replacing $X_g$ by $Y_{\sigma(g)}.$ The result now follows from \ref{groupcoset}. \endproof

\begin{lemma}\label{evlemma} Suppose $H$ is a subgroup of index $k$ in a finite group $G$ and $\rho$ is an irreducible representation of $G.$ If $\lambda\in\Lambda^{\sym}_H$ is an eigenvalue with multiplicity $m$  of $\pi^{\sym}_H(M(\rho))$ then $\lambda/k$ is an eigenvalue of $M^{\sym}(G,H)$ of multiplicity at least $m\cdot \dim(\rho).$
\end{lemma}
\demo By \fullref{rightreg} $\rho$ appears with multiplicity $\dim(\rho)$ in the right regular representation of $G.$ By \fullref{decompose} $\pi^{\sym}_H(M(G))$ has $\lambda$ as an eigenvalue of multiplicity at least $m\cdot \dim(\rho).$ The result now follows from \ref{groupcoset}.\endproof

\begin{lemma}\label{GHlemma} Suppose that $\rho\co G\rightarrow \Aut(V)$
  is a representation  and $H$ is a subgroup of $G.$ Set $P =
  \sum_{h\in H}\ \rho(h) \in \End(V).$ Then:
\begin{itemize} 
\item[\rm(i)] $\rank[P]$ equals the dimension of the subspace of $V$ on which $H$ acts trivially.
\item[\rm(ii)] $\rank\left[\pi^{\sym}_H(M(\rho))\right] = \rank[P].$
\item[\rm(iii)]
 If $\rank[P]=1$ then the non-zero eigenvalue of
 $\pi^{\sym}_H(M(\rho))$ is\break $\trace(\pi^{\sym}_H(M(\rho)))$ which is  a
 linear homogeneous polynomial in $ \Lambda^{\sym}_H$.
\end{itemize} 
\end{lemma}
\demo The image of $P$ equals the subspace of $V$ on which $H$ acts trivially which gives (i). For (ii) let $M(\rho)$ be the representation matrix of $\rho.$ Then
 $$\pi^{\sym}_H(M(\rho)) =  \sum_{g\in G}\ \rho(g) Y_{\sigma(g)}.$$
 If $a\in bH$ then $\sigma(a) = HaH\cup Ha^{-1}H = HbH\cup Hb^{-1}H=\sigma(b)$ thus $Y_{\sigma(a)}=Y_{\sigma(b)}.$ We will split the above summation up into sums over left cosets of $H.$ Let $g_1,\cdots g_k$ be a complete set of representatives of left $H$ cosets in $G$ with $g_1=1.$ Then 
 $$
 \begin{array}{rcl}
\pi^{\sym}_H(M(\rho)) & = &  \sum_{i=1}^k\sum_{h\in H}\ \rho(g_i)\rho(h) Y_{\sigma(g_i)}\\
& = &  \left(\sum_{i=1}^k \rho(g_i) Y_{\sigma(g_i)}\right) \cdot \sum_{h\in H}\ \rho(h) .
 \end{array}
 $$
We claim that the matrix $\sum_{i=1}^k \rho(g_i) Y_{\sigma(g_i)}$ is invertible.  This is because
the only  term in the sum for which $\sigma(g_i)=H$ is when $i=1.$ Thus setting $Y_H=1$ and all the other $Y's$ equal to zero gives a specialization of $\pi^{\sym}_H(M(\rho))$ which gives the identity matrix. Since this is invertible the claim, and hence part (ii) of the lemma, follows. Part (iii) follows immediately from (ii).  \endproof

\begin{proposition}\label{easyrep} Let $p$ be a prime and $G=PSL(2,{\Bbb Z}_p)$ and $H$ the subgroup of upper-triangular matrices with $1$'s on the diagonal. Then $M^{\sym}(G,H)$ has an eigenvalue, $\lambda,$ which is a non-zero rational linear combination of $\{ Y_s : s\in S \}.$ The multiplicity of $\lambda$ is at least $p.$
\end{proposition}
\demo $H$ is a cyclic subgroup of order $p$ generated by
$A = \left(
\begin{array}{cc}
1 & 1\\ 0 & 1\end{array} \right).$ There is a natural action of
$PSL(2,{\Bbb Z}_p)$ on ${\Bbb P}^1{\Bbb Z}_p.$ Let $V$ be the ${\Bbb
Q}$--vector space consisting of formal linear combinations of points in
${\Bbb P}^1{\Bbb Z}_p.$ The dimension of $V$ equals the number of
points in ${\Bbb P}^1{\Bbb Z}_p$ which is $p+1.$ Let $\rho\co PSL(2,{\Bbb
Z}_p)\rightarrow \Aut(V)$ be the natural action by
permutations. [Observe that if $g$ is any element of order $p$ in
$PSL(2,{\Bbb Z}_p)$ then $g$ must act as a $p$--cycle plus a $1$--cycle
and thus $g$ is conjugate to a power of $A.$] There is an invariant
$1$--dimensional subspace spanned by $(1,1,\cdots,1)$ on which the
action is trivial so $V = \epsilon \oplus W$ where $\epsilon$ is the
trivial $1$--dimensional representation. Since the action of $PSL(2,
\mathbb{F})$ on ${\Bbb P}^1({\Bbb Z}_p)$ is 2--transitive, $V$ splits
into two irreducible representations of dimension $1$ and $p$, by
Burrow \cite[Lemma 29.1]{Burrow}. Hence $W$ is an irreducible
representation, $\rho',$ of $G$ of dimension $p.$ The action of $H$ on
$V$ is the right regular action plus a trivial $1$--dimensional
representation. Hence the action of $H$ on $W$ is the right regular
representation. Thus the subspace of $W$ on which $H$ acts trivially
has dimension $1.$
 
 By the above calculation and \fullref{GHlemma} $\rho'$  has  an eigenvalue, $\lambda',$ which is a non-zero rational linear combination of $\{ Y_s : s\in S \}.$ Since $\rho'$ is irreducible \fullref{evlemma} implies $M^{\sym}(G,H)$ has eigenvalue $\lambda=\lambda'/p$ with multiplicity at least $\dim(\rho')=p.$\endproof

\noindent  We compute $\lambda$ for a few values of $p$:\\

\def\strutt{\vrule width 0pt depth 3pt height 14pt}
\cl{\begin{tabular}{|l||l|l|l|l|l|l|}\hline
$p$ & $\ \ \ \ \ \ \ \  \lambda$\strutt \\
\hline
$5$ &  $a +  b -  c -  d$\strutt \\
\hline
$7$ & $a + 2 b -  c -  d -  e$\strutt\\
\hline
$11$ & $a + 2 b + 2 c -  d -  e -  f -  g -  h$\strutt\\
\hline
$13$ & $a + b + 2c + 2d - e - f - g - h - i - j$\strutt\\
\hline
$17$ & $a + b + 2c + 2d + 2e - f - g - h - i - j - k - l - m$\strutt\\
\hline
$19$ & $a + 2b + 2c + 2d + 2e -f - g - h - i - j - k - l -   m - n $\strutt\\
\hline
\end{tabular}}

The following sometimes allows one to deduce information about the $m_i.$ 

\begin{proposition} Let $\rho_1,\cdots\rho_k$ be a list of the conjugacy classes of irreducible representations of $G$ and set $m_i = \rank(\pi^{\sym}_H(M(\rho_i))).$ Then $|G|/|H|\ = \sum_{\rho_i} \dim(\rho_i)\cdot m_i.$ 
\end{proposition}
\demo By \ref{rightreg} and \ref{decompose} 
$$M(G)\ =\ \bigoplus_{i=1}^k\left(\bigoplus_{j=1}^{\dim(\rho_i)}  M(\rho_i)\right)$$  Applying $\pi^{\sym}_H$ to this and taking the rank gives 
$$\rank[\pi^{\sym}_H(M(G))]\ =\ \sum_{i=1}^k \dim(\rho_i)\cdot \rank[\pi^{\sym}_H(M(\rho_i))].$$
 By \ref{grpcosetgrpdet}  the left hand side is $\rank[M^{\sym}(G,H)].$ Since $M^{\sym}(G,H)$ is non-singular  $ \rank[M^{\sym}(G,H)] =  |G|/|H|.$ \endproof

\section{Applications to 3--manifolds} \label{applications}

\proof[Proof of \fullref{maintheorem}] 
First observe that by \ref{kvhslopes} there are always slopes $\alpha,\beta$ which are not virtual homology slopes and form a basis  of $H_1(T).$  By \ref{genbdry} there is a specialization $t\co S\rightarrow{\Bbb Q}$ such that with $G,H$ as in the hypotheses the boundary matrix, $B,$  for $\tilde{Y}$ equals $ t(M^{\sym}(G,H)).$ By \ref{slopesareeigenvalues} it follows that the dimension of the $(n/m)$--eigenspace of $B$ equals the filling rank of the slope $m\alpha+n\beta$ for the cover $\tilde{Y}\rightarrow Y.$ \endproof

\begin{corollary} \label{bdrylifts} Suppose $Y$ is a compact orientable 3--manifold with  boundary a torus. Suppose that $\tilde{Y}\rightarrow Y$ is a regular cover with group of covering transformations $G$ and that $|\partial\tilde{Y}|=|G|.$ If ${\det}^{\sym}(G)$ has a linear factor of multiplicity $n$ then there is a virtual homology slope of  filling rank at least $n.$ In particular  if $G$ surjects one of ${\Bbb Z}_3,{\Bbb Z}_4,Q_8$ then there is a virtual homology slope of  filling rank at least $2.$\end{corollary}
\demo Since $|\partial\tilde{Y}|=|G|$ the stabilizer of each component of $\partial\tilde{Y}$ is trivial, thus   $H=1.$ Hence $M^{\sym}(G;H) = \pi^{\sym}(M(G))$ is the symmetrized group matrix of $G.$ The second conclusion follows from the computations in \fullref{groupdeterminants} plus the observation that if $G$ surjects $H$ then $H_1(\tilde{Y})$ surjects $H_1(\tilde{Y}_H)$ where $\tilde{Y}_H$ is the cover corresponding to $H.$ \endproof

If $\tilde{Y}\rightarrow Y$ is a regular cover with group of covering
transformations $G$ then $G$ acts on $H_1(\tilde{Y},{\Bbb Q}).$ Let
$U$ be the subspace on which the $G$--action is trivial. Then using
transfer one obtains $U\cong H_1(Y,{\Bbb Q}).$ The action of $G$ on
$V=H_1(\tilde{Y},{\Bbb Q})/U$ has no trivial summands. It follows that
$\beta_1(\tilde{Y}) = \beta_1(Y)$ or $\beta_1(\tilde{Y}) \ge
\beta_1(Y)+n$ where $n$ is the dimension of the smallest non-trivial
action of $G$ on a ${\Bbb Q}$--vector space. This gives the following
two results, which appear to have not been stated in the
literature. \fullref{closedrank} was used implicitly in Dunfield and
Thurston \cite{DT}.

\begin{theorem}\label{closedrank} Suppose that $G$ is a finite group. Let $n$ be the dimension of the smallest non-trivial ${\Bbb Q}$--representation of $G.$ Suppose that $Y$ is a  $3$--manifold and that $p\co \tilde{Y}\rightarrow Y$ is a regular cover with group of covering transformations $G.$ Either $\beta_1(\tilde{Y}) = \beta_1(Y)$ or $\beta_1(\tilde{Y}) \ge \beta_1(Y) + n.$
\end{theorem}

\noindent{\bf Example}\qua  Suppose that $G$ is simple and contains an element $g$ of  order $n.$ If $G$ acts non-trivially on a ${\Bbb Q}$ vector space $V$ then, since $G$ is simple, $g$ acts non-trivially. The irreducible ${\Bbb Q}$--representations of  $<g>$ correspond to the factorization of $x^n-1$ over ${\Bbb Q}.$ In particular if $n$ is prime then  $1+x+\cdots+x^{n-1}$ is irreducible over ${\Bbb Q}$ thus $\dim(V)\ge  n-1.$ If $p > 3$ is a prime then $PSL(2,{\Bbb Z}_p)$ is a simple group which contains an element of order $p$ so $\dim(V)\ge p-1.$  

 \begin{theorem}\label{bdryrank} Suppose that $G$ is a finite group. Let $n$ be the dimension of the smallest non-trivial ${\Bbb Q}$--representation of $G.$ Suppose that $Y$ is an orientable $3$--manifold with boundary a torus $T$ and that $p\co \tilde{Y}\rightarrow Y$ is a regular cover with group of covering transformations $G.$ Suppose that $\gamma$ is a virtual homology slope for this cover  and either $[\gamma]\ne 0\in H_1(Y)$ or $\mathrm{fillrank}(\gamma, \tilde{Y})>1.$ Then $ \mathrm{fillrank}(\gamma, \tilde{Y})\ge n.$\end{theorem}
\demo  The hypotheses on $\gamma$ imply that $\beta_1(\tilde{Y}) > \beta_1(Y).$ \endproof

\begin{lemma}\label{nonfiber} Suppose that $Y$ is a compact, connected, orientable $3$--manifold and $V$ is a subspace of $H_2(Y;{\Bbb Q})$ of dimension at least $2.$ Then there is a non-separating, connected, orientable, norm-minimizing surface $S$ in $Y$ which is not a fiber of a fibration of $Y$ over the circle and $[S] \in V.$
\end{lemma}
\demo This follows from Thurston's description of the collection of homology classes of fibers as the set of integral points in a union of cones on some of the open faces of the Thurston norm. The intersection of $V$ with such an open cone is an open cone in $V.$ There is a line of rational slope in $V$ which is in the complement of the union of these open cones. This line contains an integral point which corresponds to a norm-minimizing surface that is not a fiber.\endproof

\begin{corollary}\label{rank2} Suppose that $Y$ is a  compact, connected, orientable $3$--manifold with boundary a torus $T.$ Suppose that $p\co \tilde{Y}\rightarrow Y$ is a finite regular cover and that for each component $\tilde{T}$ of $p^{-1}T$ that the induced cover $p|\co \tilde{T}\rightarrow T$ is cyclic. Suppose that $\alpha$ is a slope on $T$ with $\mathrm{fillrank}(\alpha,\tilde{Y})\ge 2.$ Then  $Y$ has infinitely many virtually Haken Dehn-fillings.\end{corollary}
\demo Let $V$ be the subspace of $H_2(\tilde{Y};{\Bbb Q})$ generated by surfaces all of whose boundary components cover $\alpha.$ Then $\dim(V)\ge \mathrm{fillrank}(\alpha,\tilde{Y})\ge2.$ By \eqref{nonfiber} there is a norm-minimizing (hence incompressible) surface $S$ in $\tilde{Y}$ which is non-separating and not a fiber and $[S]\in V$ thus all the boundary components of $S$ cover $\alpha.$ The result follows from Lemma 7 of \cite{CooperWalsh1}.\endproof  

\proof[Proof of \fullref{virthaken} and \fullref{vhs}] 
By results of Long and Reid \cite{LongReidSimple} there is a prime $p$
such that $\pi_1Y$ surjects $G = PSL(2,{\Bbb Z}_p).$ They also show
that a surjection may be found so that $\pi_1\partial Y$ maps onto the
subgroup $H$ of upper triangular matrices with $1$'s on the
diagonal. Then by \ref{easyrep} $M^{\sym}(G,H)$ has an eigenvalue of
multiplicity at least $p$ which is a linear polynomial in
$\Lambda^{\sym}_H.$ It follows from \ref{genbdry} that the boundary
matrix for $\tilde{Y}$ has a rational eigenvalue of multiplicity at
least $p.$ This eigenvalue determines a slope $\alpha\subset T$ with
filling rank at least $p.$ This proves \fullref{vhs}.  Since $p\ge2$
\fullref{virthaken} follows from \ref{rank2}.\endproof

\begin{theorem} Suppose that $Y$ is a  compact, connected, orientable $3$--manifold with boundary a torus $T.$ Suppose that $G$ is a finite group and $\theta\co \pi_1Y\rightarrow G$ is an epimorphism and set $H=\theta(\pi_1T).$  Let $\tilde{Y}\rightarrow Y$ be the cover corresponding to $G.$ Suppose that $\det M^{\sym}(G,H)$ has a root which is a homogeneous polynomial  in $ {\Bbb Q}[\{ Y_s : s\in S \}]$ of degree one and multiplicity $m.$ Then there is a slope $\alpha$ on $T$ such that $ \mathrm{fillrank}(\alpha,\tilde{Y})\ge m.$\end{theorem} 

\begin{theorem} Suppose $Y$ is a compact, orientable 3--manifold with torus boundary and $\tilde{Y}$ is a regular cover. Then there is a basis of $H_2(\tilde{Y},\partial\tilde{Y},{\Bbb R})$ consisting of elements  whose boundaries have constant slope (possibly real rather than rational).
\end{theorem}
\demo Since the boundary matrix is symmetric there is a basis of eigenvectors.
By \ref{slopesareeigenvalues} an eigenvector gives a 2--chain with constant slope.\endproof

\subsection{An example} \label{example}
For $p>3$ the group $PSL(2, {\Bbb Z}_p)$ is simple, however $PSL(2, {\Bbb Z}_3)$ is the semi-direct product $(\mathbb{Z}_2 \oplus \mathbb{Z}_2) \rtimes \mathbb{Z}_3$. By Riley, \cite{Riley}, the fundamental group of the figure-8 knot complement surjects  $PSL(2, {\Bbb Z}_3)$ with the peripheral subgroup mapping onto $H$. Now we claim that there is only one $(\mathbb{Z}_2 \oplus \mathbb{Z}_2) \rtimes \mathbb{Z}_3$ cover of the figure-8 knot complement such that the peripheral group maps to $Z_3$.  Let $\Gamma$ be the fundamental group of the figure-8 knot complement.  Then a surjection $\Gamma \rightarrow  (\mathbb{Z}_2 \oplus \mathbb{Z}_2) \rtimes \mathbb{Z}_3$ induces a surjection  $\Gamma \rightarrow \mathbb{Z}_3$.  This map factors through the abelianization, so it corresponds to the unique 3--fold cyclic cover of the figure-8 knot complement. The first homology of this manifold is $\mathbb{Z} \oplus \mathbb{Z}_4 \oplus \mathbb{Z}_4$.  A map from the fundamental group of this manifold to $\mathbb{Z}_2 \oplus \mathbb{Z}_2$ factors through the abelianization, and is unique when the boundary torus lifts to the associated cover $\tilde Y$.   Any filling of the figure-8 knot complement that lifts to the three-fold cyclic cover will lift to $\tilde Y$. The rational homology of $\tilde Y$ is carried by the boundary. 

According to \ref{PSL2cover}, there is a virtual homology slope of
filling rank at least 3 for the associated regular cover.  We claim
that this slope is a meridian.  Indeed, orbifold-filling along a (3,0)
curve results in a Euclidean orbifold (Cooper, Hodgson and Kerckhoff
\cite{Orbifoldbook}). The unique three-fold cyclic manifold cover of
this is a orientable Euclidean manifold $Y$ with $H_1(Y) =
\mathbb{Z}_4 \oplus \mathbb{Z}_4$.  By the classification of such
manifolds in Wolf \cite{Wolfconstant}, $\pi_1(Y)$ is the group $G_6$,
and this manifold has a unique $\mathbb{Z}_2 \oplus \mathbb{Z}_2$
cover which is the 3--torus.  This is an equivariant Dehn-filling of
$\tilde Y$, where the filling slopes are pre-images of the meridian of
the complement of the figure-8 knot.  Therefore, this is a virtual
homology slope of filling rank 3.

\bibliographystyle{gtart}
\bibliography{link}

\end{document}